\documentclass[12pt]{article}
\usepackage{amssymb,amsfonts,amsmath,psfrag,eepic,colordvi,graphicx,epsfig}
\usepackage{amssymb,latexsym,graphics,array}
\usepackage{MnSymbol}
\usepackage[enableskew]{youngtab}
\usepackage{float,enumerate,enumitem}
\usepackage{tikz,ifthen}
\usepackage{pgfplots}
\usetikzlibrary{math}

\usepackage{framed}
\usepackage{hyperref}
\hypersetup{colorlinks=true}

\topskip  
\newcommand{\rmnum}[1]{\romannumeral #1}

\numberwithin{equation}{section}
\newtheorem{theorem}{Theorem}[section]

\newtheorem{conjecture}[theorem]{Conjecture}
\newtheorem{remark}[theorem]{Remark}
\newtheorem{lemma}[theorem]{Lemma}

\begin{document}
	\parskip 6pt
	
	\pagenumbering{arabic}
	\def\sof{\hfill\rule{2mm}{2mm}}
	\def\ls{\leq}
	\def\gs{\geq}
	\def\SS{\mathcal S}
	\def\qq{{\bold q}}
	\def\MM{\mathcal M}
	\def\TT{\mathcal T}
	\def\EE{\mathcal E}
	\def\lsp{\mbox{lsp}}
	\def\rsp{\mbox{rsp}}
	\def\pf{\noindent {\it Proof.} }
	\def\mp{\mbox{pyramid}}
	\def\mb{\mbox{block}}
	\def\mc{\mbox{cross}}
	\def\qed{\hfill \rule{4pt}{7pt}}
	\def\block{\hfill \rule{5pt}{5pt}}
	\def\lr#1{\multicolumn{1}{|@{\hspace{.6ex}}c@{\hspace{.6ex}}|}{\raisebox{-.3ex}{$#1$}}}
	\def\red{\textcolor{red}}

\begin{center}
	{\Large\bf  On a conjecture concerning the $r$-Euler-Mahonian statistic on permutations}
\end{center}

\begin{center}
 
	{\small Kaimei Huang$^a$,  Zhicong Lin$^b$,
		Sherry H.F. Yan$^{a,*}$\footnote{$^*$Corresponding author.}\footnote{{\em E-mail address:} hfy@zjnu.cn. }}

 	$^{a}$Department of Mathematics,
 Zhejiang Normal University\\
 Jinhua 321004, P.R. China
 
 	$^{b}$ Research Center for Mathematics and Interdisciplinary Sciences,  Shandong University \& Frontiers Science Center for Nonlinear Expectations, Ministry of Education, Qingdao 266237, P.R. China

\end{center}

\noindent {\bf Abstract.} A pair  $(\mathrm{st_1}, \mathrm{st_2})$  of permutation statistics is said to be   $r$-Euler-Mahonian  if  $(\mathrm{st_1}, \mathrm{st_2})$ and 
  $( \mathrm{rdes}$, $\mathrm{rmaj})$ are equidistributed over the set  $\mathfrak{S}_{n}$ of all permutations of $\{1,2,\ldots, n\}$, where $\mathrm{rdes}$ denotes the $r$-descent number  and $\mathrm{rmaj}$ denotes the $r$-major index introduced by Rawlings. The main objective of this paper is to   prove that  $(\mathrm{exc}_r, \mathrm{den}_r)$  and  $( \mathrm{rdes}$, $\mathrm{rmaj})$ are equidistributed over  $\mathfrak{S}_{n}$, thereby confirming a recent conjecture posed by Liu. When $r=1$, the result recovers the equidistribution of $(\mathrm{des}, \mathrm{maj})$ and $(\mathrm{exc}, \mathrm{den})$, which was first conjectured by Denert and proved by Foata and Zeilberger. 




\section{Introduction}
   Let $\mathfrak{S}_{n}$ denote the set of permutations of $[n]:=\{1,2,\ldots, n\}$. A permutation $\pi\in \mathfrak{S}_{n}$ is usually written in one-line notation as $\pi=\pi_1\pi_2\cdots \pi_n$. Given a permutation  $\pi\in \mathfrak{S}_{n}$,  an index $i$, $1\leq i< n$,  is called a  {\em descent  } of $\pi$ if $\pi_i>\pi_{i+1}$.    
 Define  the {\em descent set}  of $\pi$  to be
$$
\mathrm{Des}(\pi)=\{i\in[n-1]\mid \pi_i>\pi_{i+1}\}.
$$
 Let $\mathrm{des(\pi)}=|\mathrm{Des}(\pi)|$. 
  The {\em major index} of $\pi$, denoted by $\mathrm{maj(\pi)}$, is defined to be
$$
\mathrm{maj(\pi)}=\sum\limits_{i\in \mathrm{Des(\pi)} } i.
$$
For example, if we let  $\pi=475398216$, then $\mathrm{Des}(\pi)=\{2,3,5,6,7\}$, $\mathrm{des}(\pi)=5$ and $\mathrm{maj}(\pi)=2+3+5+6+7=23$.
An index $i$, $1\leq i< n$, is called an {\em excedance place}
 of $\pi$ if $\pi_i>i$. The {\em excedance number} of $\pi$, denoted by $\mathrm{exc}(\pi)$,  is defined to be the number of excedance places of $\pi$.  For example, if we let  $\pi=475398216$, then the excedance places of $\pi$ are $1,2,3,5,6$  and $\mathrm{exc}(\pi)=5$.  The classical {\em Eulerian polynomial} (see~\cite[pp.~32-33]{ST})  $A_n(t)$ is  defined as 
 $$
 A_n(t)=\sum_{\pi\in \mathfrak{S}_{n}}t^{\mathrm{des}(\pi)}=\sum_{\pi\in \mathfrak{S}_{n}}t^{\mathrm{exc}(\pi)}.
 $$
Given a permutation  $\pi \in \mathfrak{S}_{n}$, a pair $(i,j)$ is called an {\em inversion pair} of $\pi$ if $i<j$ and $\pi_i>\pi_j$.  Let $\mathrm{Inv}(\pi)$ denote the set of all inversion pairs of $\pi$.
The {\em inversion number} of $\pi$, denoted by $\mathrm{inv(\pi)}$, is defined to be the number of inversion pairs  of $\pi$.   
 MacMahon's equidistribution  theorem \cite{Mac} states that
\begin{equation}\label{invmaj}
\sum\limits_{\pi\in \mathfrak{S}_{n}}  q^\mathrm{inv(\pi)}=\sum\limits_{\pi\in  \mathfrak{S}_{n }}q^\mathrm{maj(\pi)}=[n]_q[n-1]_q\ldots [1]_q,
\end{equation}
where $[k]_q:=1+q+q^2+\ldots+q^{k-1}$. 
The problem of finding MacMahon type results  for other combinatorial objects has been extensively  investigated in the literature. For example, see~\cite{Chen1} for matchings and set partitions, \cite{Chen2} for $01$-fillings of moon polyominoes, \cite{Liu1, Liu2, Yan} for $k$-Stirling permutations,  and~\cite{Wilson} for ordered set partitions.

A pair  $(\mathrm{st_1}, \mathrm{st_2})$  of permutation statistics is said to be {\em Euler-Mahonian} if 
$$
\sum\limits_{\pi\in \mathfrak{S}_{n}}  t^{\mathrm{st_1}(\pi)}q^{\mathrm{st_2}(\pi)}=  \sum\limits_{\pi\in \mathfrak{S}_{n}}  t^{\mathrm{des}(\pi)}q^{\mathrm{maj}(\pi)}.
$$
 Denert \cite{Denert} introduced a new permutation statistic   $\mathrm{den}$ and  posed a conjecture which asserts that  $\mathrm{(exc, den)}$ is Euler-Mahonian, that is,
 \begin{equation}\label{den-maj}
 	 \sum_{\pi\in \mathfrak{S}_{n}}t^{\mathrm{des}(\pi)}q^{\mathrm{maj}(\pi)}=\sum_{\pi\in \mathfrak{S}_{n}}t^{\mathrm{exc}(\pi)}q^{\mathrm{den}(\pi)}.
 \end{equation}
  This conjecture was first  proved by Foata and Zeilberger \cite{Foata-Zeilberger}  and  a bijective proof was provided by Han \cite{Han}. Here we adopt an equivalent definition of Denert's statistic
 due to Foata and Zeilberger \cite{Foata-Zeilberger}.  Given a permutation  $\pi\in \mathfrak{S}_{n}$, if  $\pi_i>i$, then  $\pi_i$ is called an {\em excedance letter} of $\pi$. Let $\mathrm{Excp}(\pi)$ and $\mathrm{Nexcp}(\pi)$ denote the set of   excedance places of $\pi$ and the set of non-excedance places of $\pi$, respectively.   Assume that  $\mathrm{Excp}(\pi)=\{i_1, i_2, \ldots, i_k\}$ with $i_1<i_2<\ldots<i_k$ and   $\mathrm{Nexcp}(\pi)=\{j_1, j_2, \ldots, j_{n-k}\}$
 with $j_1<j_2<\ldots<j_{n-k}$.  Define
 $ 
 \mathrm{Exc}(\pi)=\pi_{i_1}\pi_{i_2}\ldots \pi_{i_k} 
 $ 
 and 
 $ 
 \mathrm{Nexc}(\pi)=\pi_{j_1}\pi_{j_2}\ldots \pi_{j_{n-k}}.
 $ 
  Denert's statistic of $\pi$, denoted by $\mathrm{den}(\pi)$, is defined as
  $$
  \mathrm{den}(\pi)=\sum\limits_{i\in \mathrm{Excp}(\pi)}i+{\mathrm{inv}}(\mathrm{Exc}(\pi))+{\mathrm{inv}}(\mathrm{Nexc}(\pi)).
  $$
  For example, if we let  $\pi=475398216$, then  $\mathrm{Excp}(\pi)=\{1,2,3,5,6\}$, 
  $\mathrm{Exc}(\pi)=47598$,  and $\mathrm{Nexc}(\pi)=3216$.
  Then $\mathrm{den}(\pi)$ is given by 
  $$
  \begin{array}{ll}
  	\mathrm{den}(\pi)&=\sum\limits_{i\in \mathrm{Excp}(\pi)}i+{\mathrm{inv}}(\mathrm{Exc}(\pi))+{\mathrm{inv}}(\mathrm{Nexc}(\pi))\\
     &=1+2+3+5+6+{\mathrm{inv}}(47598)+{\mathrm{inv}}(3216)\\
    
  &=22.
  \end{array}
  $$

      Let $r\geq 1$. Given a permutation $\pi\in \mathfrak{S}_{n}$, let 
$$
\mathrm {rDes}(\pi)=\{i\in [n-1]\mid \pi_{i}\geq \pi_{i+1}+r\},
$$
$$
\mathrm {rInv}(\pi)=\{(i,j)\in \mathrm{Inv}(\pi)\mid  \pi_{i}<\pi_{j}+r  \},
$$
An index $i$ is said to be an {\em $r$-descent}  if $i\in \mathrm {rDes}(\pi) $. The {\em $r$-descent  number}, denoted by $\mathrm{rdes}(\pi)$, is defined to be the number  of $r$-descents.
The {\em   $r$-major
	index} of $\pi$, denoted by $\mathrm{rmaj}(\pi)$, is defined to be 
$$
\mathrm{rmaj}(\pi)=\sum\limits_{i\in \mathrm{rDes}(\pi)}i +|\mathrm{rInv}(\pi)|. 
$$
For example, if we let    $\pi=475398216$ and $r=2$, then $\mathrm{rDes}(\pi)=\{2, 3, 6\}$, $\mathrm{rdes}(\pi)=3$ and $\mathrm{rInv}(\pi)=\{(1,4), (2,9), (4, 7), (5,6),$ $ (7,8)\}$. Thus,  $\mathrm{rmaj}(\pi)$ is given by
$$
\begin{array}{ll}
	\mathrm{rmaj}(\pi)  &=\sum\limits_{i\in \mathrm{rDes}(\pi)}i +|\mathrm{rInv}(\pi)|\\
	&=2+3+6+5\\
	&=16.
\end{array}
$$
The notions of $r$-descent number and   $r$-major
index were introduced by
Rawlings~\cite{Raw}. Clearly, $\mathrm{rmaj}$ reduces to $\mathrm{maj}$ when $r=1$, and reduces to $\mathrm{inv}$ when $r\geq n$. Thus the statistic     $\mathrm{rmaj}$ interpolates between  $\mathrm{maj}$ and $\mathrm{inv}$.
Rawlings~\cite{Raw}  proved bijectively that $\mathrm{rmaj}$ and $\mathrm{inv}$ are equidistributed over $\mathfrak{S}_{n}$.
A pair  $(\mathrm{st_1}, \mathrm{st_2})$  of permutation statistics is said to be {\em $r$-Euler-Mahonian} if 
$$
\sum\limits_{\pi\in \mathfrak{S}_{n}}  t^{\mathrm{st_1}(\pi)}q^{\mathrm{st_2}(\pi)}=  \sum\limits_{\pi\in \mathfrak{S}_{n}}  t^{\mathrm{rdes}(\pi)}q^{\mathrm{rmaj}(\pi)}.
$$
  
   Recently,  Liu~\cite{Liujcta} introduced the notions of {\em $r$-level excedance number} and  {\em $r$-level Denert's statistic}. Fix an integer $r\geq 1$.  Given a permutation $\pi\in \mathfrak{S}_{n}$, an index $i\in[n-1]$ is called an {\em $r$-level excedance place}  and $\pi_i$ is called an {\em $r$-level excedance letter} if $\pi_i>i$ and $\pi_i\geq r$. For example, if we let    $\pi=475398216$ and $r=5$, then the $5$-level excedance places of $\pi$ are $2,3, 5,6$ and $5$-level excedance letters of $\pi$ are $5,7,8,9$. 
   Let $\mathrm{Excp}_r(\pi)$ and $\mathrm{Nexcp}_r(\pi)$ denote the set of $r$-level excedance places of $\pi$ and the set of non-$r$-level excedance places of $\pi$, respectively.  Let $\mathrm{Excl}_r(\pi)$ denote the set of $r$-level excedance letters of $\pi$. 
     Assume that  $\mathrm{Excp}_r(\pi)=\{i_1, i_2, \ldots, i_k\}$ with $i_1<i_2<\ldots<i_k$ and   $\mathrm{Nexcp}_r(\pi)=\{j_1, j_2, \ldots, j_{n-k}\}$
       with $j_1<j_2<\ldots<j_{n-k}$.  Define
   $ 
   \mathrm{Exc}_{r}(\pi)=\pi_{i_1}\pi_{i_2}\ldots \pi_{i_k},
   $ 
   and 
   $ 
   \mathrm{Nexc}_{r}(\pi)=\pi_{j_1}\pi_{j_2}\ldots \pi_{j_{n-k}}.
   $ 
    The {\em  $r$-level Denert's statistic} of $\pi$, denoted by $\mathrm{den}_r(\pi)$, is defined as
   $$
   \mathrm{den}_{r}(\pi)=\sum\limits_{i\in \mathrm{Excp}_r(\pi)} i+{\mathrm{inv}}(\mathrm{Exc}_{r}(\pi))+{\mathrm{inv}}(\mathrm{Nexc}_{r}(\pi)).
   $$
    For example, if we let    $\pi=475398216$ and $r=5$, then $\mathrm{Excp}_5(\pi)=\{2,3,5,6\}$,   $\mathrm{Excl}_5(\pi)=\{5,7,8,9\}$, $\mathrm{Nexcp}_5(\pi)=\{1,4,7,8,9\}$, 
   $\mathrm{Exc}_{5}(\pi) =7598$, and  $\mathrm{Nexc}_{5}(\pi) =43216$. Therefore,   $\mathrm{den}_5(\pi)$ is  given by 
   $$
   \begin{array}{ll}
   	\mathrm{den}_5(\pi)&=\sum\limits_{i\in \mathrm{Excp}_5(\pi)} i+{\mathrm{inv}}(\mathrm{Exc}_{5}(\pi))+{\mathrm{inv}}(\mathrm{Nexc}_{5}(\pi))\\
   	&= 2+3+5+6+{\mathrm{inv}}( 7598)+{\mathrm{inv}}(43216)\\
      	&=24.
   \end{array}
   $$   
    It should be mentioned that  $\mathrm{den}_r$ reduces to  inv when $r\geq n$,  and  reduces to den when $r=1$.   Thus the statistic     $\mathrm{den}_r$ interpolates between   $\mathrm{inv}$ and $\mathrm{den}$. 
  The {\em $r$-level excedance number } of $\pi$,  denoted by $\mathrm{exc}_r(\pi)$,  is defined to be the number of $r$-level excedance letters that occurs weakly to the right of $\pi_{r}$, that is, 
    $$\mathrm{exc}_r(\pi)=|\{i\mid \pi_i>i\geq r\}|.$$  For example, if we let    $\pi=475398216$ and $r=5$, then 
    $\mathrm{exc}_{5}(\pi) =2$. 
   Since $\{i\mid \pi_i>i\geq r\} \subseteq \{i\mid \pi_i>i, \pi_i\geq r\}=\mathrm{Excp}_r(\pi)$ and $|\mathrm{Exc}_{r}(\pi)|=|\mathrm{Excp}_r(\pi)|$,  we have 
   $\mathrm{exc}_r(\pi) \leq |\mathrm{Exc}_{r}(\pi)|$, where
   $|\mathrm{Exc}_{r}(\pi)|$ denotes the length of the word $\mathrm{Exc}_{r}(\pi)$.  Clearly,  the statistic $\mathrm{exc}_r$ reduces to exc when $r=1$.

   Motivated by  finding  a generalization of~\eqref{den-maj} for $r$-Euler-Mahonian statistics, Liu~\cite{Liujcta} posed the following conjecture.   
   \begin{conjecture}[Liu~\cite{Liujcta}]\label{con1}     
   	For all $r\geq 1$, the pair $(\mathrm{exc}_r, \mathrm{den}_r)$ is $r$-Euler-Mahonian, that is, 
   	\begin{equation}\label{rden-rmaj}
   		\sum_{\pi\in \mathfrak{S}_{n}}t^{\mathrm{rdes}(\pi)}q^{\mathrm{rmaj}(\pi)}=\sum_{\pi\in \mathfrak{S}_{n}}t^{\mathrm{exc}_r(\pi)}q^{\mathrm{den}_r(\pi)}.
   	\end{equation}
   	\end{conjecture}
   The main objective of this paper is to prove that $(\mathrm{rdes}, \mathrm{rmaj})$ and $(\mathrm{exc}_r, \mathrm{den}_r)$ are equidistributed over $\mathfrak{S}_{n}$, thereby confirming Conjecture \ref{con1}. 
   \section{Two insertion methods}
    
   \subsection{Rawlings's insertion method}
   In this subsection, we give a review of Rawlings's insertion method \cite{Raw}. Let $r\geq 1$. Given $\pi\in \mathfrak{S}_{n-1}$, let 
   $$
   S(\pi)=\{\pi_{i}\mid \pi_{i-1}\geq \pi_i+r\}\cup \{\pi_{i}\mid \pi_i>n-r\}. 
   $$
   Clearly, we have $|S(\pi)|=\mathrm{rdes}(\pi)+r-1$. The {\em $r$-maj-labeling} of $\pi$  is obtained as follows.
   \begin{itemize}
   	\item Label the  space after $\pi_{n-1}$ by $0$.
   	\item Label the spaces before those letters of $\pi$ that are belonging to $S(\pi)$ from right to left with $1,2,\ldots, \mathrm{rdes}(\pi)+r-1$. 
   	\item Label the remaining spaces from left to right with $\mathrm{rdes}(\pi)+r, \ldots, n-1$.
   \end{itemize}

For example, if we   let    $\pi=475398216$ and $r=2$, then we have $\mathrm{rDes}(\pi)=\{2, 3, 6\}$ and  $S(\pi)=\{2, 3,5, 9 \}$,  and thus  the $2$-maj-labeling of $\pi$ is given by
$$
\begin{array} {llllllllll}
_{5}4&_67&_4 5&_3 3&_29&_78&_12&_81&_96&_0,
\end{array}
$$
where the labels of the spaces are written as subscripts. 
Define $$\phi^{\mathrm{maj}}_{r,n}:\mathfrak{S}_{n-1}\times  \{0, 1, \ldots, n-1\} \longrightarrow \mathfrak{S}_{n}$$ by mapping the pair $( \pi, c)$ to the  permutation in $\mathfrak{S}_{n}$  obtained by inserting   $ n$ at the space in $\pi$ which  is labeled by $c$ in the $r$-$\mathrm{maj}$-labeling of $\pi$. Continuing with our running example, we have
$\phi^{\mathrm{maj}}_{2,10}(475398216, 7)=47539(10)8216$.

Rawlings \cite{Raw} proved that the map $\phi^{\mathrm{maj}}_{r,n}$  verifies  the following celebrated  properties.

\begin{lemma}[Rawlings~\cite{Raw}] \label{Raw}
	The map $\phi^{\mathrm{maj}}_{r,n}:\mathfrak{S}_{n-1}\times  \{0, 1, \ldots, n-1\} \longrightarrow \mathfrak{S}_{n}$ is a bijection  such that for $\pi\in \mathfrak{S}_{n-1}$ and $0\leq c\leq n-1$, we have 
	$$
	\mathrm{rdes}(\phi^{\mathrm{maj}}_{r,n}(\pi, c))=\left\{ \begin{array}{ll}
		\mathrm{rdes}(\pi)&\, \mathrm{if}\,\, 0\leq c\leq \mathrm{rdes}(\pi)+r-1,\\
			\mathrm{rdes}(\pi)+1&\, \mathrm{otherwise}, 
		\end{array}
	\right.
	$$
	and 
	$$
	\mathrm{rmaj}(\phi^{\mathrm{maj}}_{r,n}(\pi, c))=\mathrm{rmaj}(\pi)+c.
	$$
\end{lemma}
When $r=1$, the map $\phi^{\mathrm{maj}}_{r,n}$ coincides with Carlitz's insertion method  \cite{Carlitz}  which was employed  to prove  the equidistribution of $ \mathrm{maj} $ and $ \mathrm{inv}$ over $\mathfrak{S}_{n}$. 

\subsection{Han's insertion method}
In the subsection, we give a description of Han's insertion method introduced in~\cite{Han}. Here we adopt Liu's labeling way \cite{Liujcta} to  describe  Han's insertion method. Given $\pi\in \mathfrak{S}_{n-1}$, assume that $\mathrm{exc}(\pi)=s$ and that $e_1<e_2<\cdots<e_s$ are the excedance letters of $\pi$. The {\em den-labeling of the index sequence} of $\pi$ is obtained by labelling the spaces of the index sequence $12 \cdots (n-1)$ of $\pi$ as follows.
 \begin{itemize}
	\item Label the  space after the index $n-1$ by $0$.
	\item Label the space  before the index $e_i $ with $s+1-i$ for all $1\leq i\leq s$.  
	\item Label the remaining spaces from left to right with $s+1, \ldots, n-1$.
\end{itemize}
The {\em  insertion letter} of the index sequence of $\pi$ is obtained by labeling  the space weakly to the left of the index $e_1$ with $e_1$, labeling  the space weakly to the left of the index $e_2$ and to the right of the index $e_1$ with $e_2$, $\ldots$, and labeling the space to the right of the index $e_s$ with $n$. 

For example, let    $\pi=836295417$. Then the  excedance letters are $e_1=3$, $e_2=6$, $e_3=8$ and  $e_4=9$ and thus the den-labeling  and the  insertion letters    of the index sequence of  $\pi$ are given by$$
\begin{array}{lllllllllll} 
	&_{5}1&_{6}2&_{4}{\bold 3}&_{7}4&_{8}5&_{3}{\bold 6}&_{9}7&_{2}{\bold 8}&_{1}{\bold 9}&_{0}\\
	& \uparrow&\uparrow &\uparrow &\uparrow &\uparrow&\uparrow&\uparrow&\uparrow&\uparrow&\uparrow\\
	&3 &3 &3 &6 &6&6&8&8&9&10,
\end{array}
$$
where the labels of the spaces in the den-labeling are written as subscripts and the insertion letters of the spaces  are written on the bottom row.

Define $$\phi_n: \mathfrak{S}_{n-1}\times \{0, 1, \ldots, n-1\} \longrightarrow \mathfrak{S}_{n}$$ by mapping the pair $(  \pi, c)$ to the  permutation in $\mathfrak{S}_{n}$  according to the following rules.
\begin{itemize}
 \item If $c=0$, then let $\phi_n(\pi,c)$ be the permutation obtained from $\pi$ by inserting $n$ immediately  after $\pi_{n-1}$. 
 \item For $c\geq 1$, assume that the space labeled by $c$ in the den-labeling of the index sequence of $\pi$ is immediately before the index $x$, and assume that the insertion letter for this space is $e_{y}$, where $1\leq y\leq \mathrm{exc}(\pi)+1$ and $e_{\mathrm{exc}(\pi)+1}=n$. Define $\pi'=\pi'_1\pi'_2\ldots \pi'_{n-1}$ from $\pi$ by replacing $e_j$ with  $e_{j+1}$ for all $y\leq j\leq \mathrm{exc}(\pi) $. Then define $\phi_n(\pi,c)$ to be the permutation obtained from $\pi'$ by inserting the letter $e_y$ immediately before the letter $\pi'_{x}$.  
\end{itemize}

Take  $\pi=836295417$ and $c=8$. Clearly, the space labeled by $8$ in the den-labeling of the index sequence of $\pi$ is immediately before the index $x=5$, and   the insertion letter for this space is $e_2=6$.  Then we get a permutation $\pi'= 9382(10)5417$ from $\pi$ by replacing replacing $e_2=6$ with $e_3=8$,  replacing $e_3=8$ with $e_4=9$, and replacing $e_4=9$ with $e_5=10$. Finally, we get $\phi_{10}(\pi,8)=93826(10)5417 $  from $\pi'$ by inserting $e_2=6$ immediately before $\pi'_x=\pi'_5=10$.


Following Han~\cite{Han}, an index $i\in[n]$ is called a {\em grande-fixed place} of $\tau\in\mathfrak{S}_n$ if 
\begin{itemize}
\item either $\tau_i=i$ (i.e., $i$ is a fixed place)
\item  or $\tau_i>i$ and $\{\tau_j \mid \tau_j>j\text{ and }i\leq \tau_j<\tau_i\}=\emptyset$.
\end{itemize}
For example, if $\tau=816259437$, then the indices   $3$ and $5$ are grande-fixed places, while the indices $1$ and   $6$ are not grande-fixed places.

Han~\cite{Han} proved that the map $\phi_{n}$ possesses  the following desired properties.

\begin{lemma}[Han~\cite{Han}] \label{Han}	 
	The map $\phi_{n}: \mathfrak{S}_{n-1}\times \{0, 1, \ldots, n-1\} \longrightarrow \mathfrak{S}_{n}$ is a bijection  such that 
	for $\pi\in \mathfrak{S}_{n-1}$ and $0\leq c\leq n-1$,  the resulting permutation $\tau=\phi_n(\pi)$ verifies the following properties.
	\begin{itemize}
		\item[\upshape(\rmnum{1})]  We have
		$$ 
		\mathrm{exc}(\tau)=\left\{ \begin{array}{ll}
			\mathrm{exc}(\pi)&\, \mathrm{if}\,\, 0\leq c\leq \mathrm{exc}(\pi),\\
			\mathrm{exc}(\pi)+1&\, \mathrm{otherwise}, 
		\end{array}
		\right.
		$$
		and 
		$$
		\mathrm{den}(\tau)=\mathrm{den}(\pi)+c.
		$$

		\item[\upshape(\rmnum{2})] Assume that the space before the index $x$ is labeled by $c$ under   the den-labeling of the index sequence of $\pi$  when $c>0$.  Then, the index $x$  is the greatest  grande-fixed place   of $\tau$ when $c>0$ and  the index $n$  is the greatest      grande-fixed place   of $\tau$ when $c=0$.  
	\end{itemize}
\end{lemma}
Relying on Lemma \ref{Han},  Han \cite{Han} provided a bijective proof of the equidistribution of  $(\mathrm{exc}, \mathrm{den})$ and $(\mathrm{des}, \mathrm{maj})$.

\section{Proof of Conjecture \ref{con1}}

This section is devoted to the proof of Conjecture \ref{con1}. Our proof is based on introducing a  nontrivial extension of Han's insertion method. 

Let $1\leq r< n$. 
Given $\pi\in \mathfrak{S}_{n-1}$, assume that $|\mathrm{Exc}_r(\pi)|=s$ and that $\mathrm{Excl}_r(\pi)=\{e_1, e_2, \ldots, e_s\}$ with $e_1<e_2<\ldots<e_s$. 
Define 
$$
\mathcal{A}(\pi)=\{i\mid   \pi_i\geq r, i<r\}
\quad\text{and}\quad 
\mathcal{B}(\pi)=\{i\mid   \pi_i< r, i<r\}.
$$
Notice that 
$|\mathcal{A}(\pi)|+\mathrm{exc}_r(\pi)=s$ and $|\mathcal{A}(\pi)|+|\mathcal{B}(\pi)|=r-1$. Thus, we have $|\mathcal{B}(\pi)|+s=\mathrm{exc}_r(\pi)+r-1$.  Clearly, we have $\mathcal{B}(\pi)\cap \{e_1, e_2, \ldots, e_s\}=\emptyset $ as $e_i\geq r$ for all $1\leq i\leq s$.

 The {\em $r$-level-den-labeling} of    the index sequence  of $\pi$ is obtained by labeling the spaces of the index sequence $12\ldots (n-1)$ of $\pi$ as follows.
\begin{itemize}
	\item Label the  space after the index $n-1$ by $0$.
	\item Label the space  before the index $e_i $ with $s+1-i$ for all $1\leq i\leq s$.  
		\item Label the space  before the index of $\pi$ belonging to $\mathcal{B}(\pi)$ from left to  right with $s+1, \ldots, \mathrm{exc}_r(\pi)+r-1$.
	\item Label the remaining spaces from left to right with $\mathrm{exc}_r(\pi)+r, \ldots, n-1$.
\end{itemize}

For example, let    $\pi=836295417$ and $r=5$. Then the $5$-level excedance letters are $e_1=6$, $e_2=8$ and  $e_3=9$. Clearly, we have $\mathcal{A}(\pi)=\{1,3\}$ and $\mathcal{B}(\pi)=\{2,4\}$. Then the $5$-level-den-labeling    of the index sequence of  $\pi$ is given by$$
\begin{array}{lllllllllll} 
	&_{6}1&_{4}2&_{7}3&_{5}4&_{8}5&_{3}{\bold 6}&_{9}7&_{2}{\bold 8}&_{1}{\bold 9}&_{0}, 
		 
\end{array}
$$
where the labels of the spaces in the $5$-level-den-labeling are written as subscripts.

 \begin{remark}\label{ob1}
	For  $x\geq r$, the label of the space before the index $x$ has the same label in the $r$-level-den-labeling of   the index sequence of $\pi$  as that in den-labeling of   the index sequence of $\pi$.
\end{remark}

\begin{theorem}\label{thden}
	Let $1\leq r< n$.
 There is a bijection
  $$\phi^{\mathrm{den}}_{r,n}:  \mathfrak{S}_{n-1}\times \{0, 1, \ldots, n-1\} \longrightarrow \mathfrak{S}_{n}$$ such that  for $\pi\in \mathfrak{S}_{n-1}$ and $0\leq c\leq n-1$,  we have 
  
  \begin{equation}\label{eqprop1}
  \mathrm{exc}_r(\phi^{\mathrm{den}}_{r,n}( \pi, c))=\left\{ \begin{array}{ll}
  	\mathrm{exc}_r(\pi)&\, \mathrm{if}\,\, 0\leq c\leq \mathrm{exc}_r(\pi)+r-1,\\
  	\mathrm{exc}_r(\pi)+1&\, \mathrm{otherwise}, 
  \end{array}
  \right.
  \end{equation}
  and 
  \begin{equation}\label{eqprop2}
  \mathrm{den}_r(\phi^{\mathrm{den}}_{r,n}( \pi, c))=\mathrm{den}_r(\pi)+c.
  \end{equation} 
\end{theorem}

Before we establish the bijection  $\phi^{\mathrm{den}}_{r,n}$, we shall define two maps   which will play
essential roles in the construction of $\phi^{\mathrm{den}}_{r,n}$. 

 Let $1\leq r<n$. An index $i$ is called an {\em $r$-level grande-fixed place} of $\tau\in\mathfrak{S}_{n}$ if $i\geq r$ and $i$ is a grande-fixed place of $\tau$. For any positive integers $a$ and $b$ with $a< b$,  let $(a,b)$ denote the set of all integers $i$ with $a<i<b$  and let $[a,b):=(a,b)\cup\{a\}$. Then, the index $i$ with $i\geq r$ is an $r$-level grande-fixed place of $\tau$ if and only if  either $\tau_i=i$  or $\tau_i>i$ and $[i, \tau_i)\cap \mathrm{Excl}_r(\tau)=\emptyset$.  For example, if $\tau=436259817$ and $r=3$, then we have $\mathrm{Excl}_3(\tau)=\{3, 4,6,8,9\}$. It is easily seen that  the indices $5$ and  $7$ are $3$-level grande-fixed places of $\tau$, while the indices $3$ and  $6$ are  not  $3$-level grande-fixed places  of $\tau$. 
 
 Let $\mathfrak{S}_{n, r}$ be the subset of permutations in $\mathfrak{S}_{n}$ without any  $r$-level grande-fixed place.   Define
 $$
 \mathfrak{S}^*_{n, r}=\{\tau\in \mathfrak{S}_{n,r}\mid r\in \mathrm{Excp}_r(\tau)\}.
 $$
Let $\mathcal{X}_{n,r}$ (resp., $\mathcal{Y}_{n,r}$) be the subset of pairs $(\pi, c)\in \mathfrak{S}_{n-1}\times \{ 1, \ldots,  n-1\}$ such that the space labeled by $c$  under the $r$-level-den-labeling of the index sequence of $\pi$    is immediately  before the index $x$   for some $ x\in\mathcal{A}(\pi)$ (resp.,~for some $ x\in\mathcal{B}(\pi)$).

 \begin{framed}
 \begin{center}
{\bf The map $\alpha_{r,n}:\mathcal{X}_{n,r}\longrightarrow \mathfrak{S}^*_{n, r}$} 
\end{center}
 Let $(\pi, c)\in \mathcal{X}_{n,r}$ with   $|\mathrm{Exc}_r(\pi)|=s$. Assume that   the space labeled by $c$  under the $r$-level-den-labeling of the index sequence of $\pi$    is  immediately before the index $x$   and that $\mathrm{Excl}_r(\pi)=\{e_1,e_2, \ldots, e_{s}\}$ with $e_1<e_2\ldots<e_s$.   
Suppose that $\mathcal{A}(\pi)=\{i_1, i_2, \ldots, i_k\}$ with $i_1<i_2<\ldots<i_k$ and that $x=i_{y}$ for some $1\leq y\leq k$. 
Define  $\alpha_{r,n}(\pi,c)$ to be the permutation  constructed by the following procedure.
\begin{itemize}
	\item Step 1: Construct a permutation  $\sigma=\sigma_1\sigma_2\ldots \sigma_{n-1}$  from $\pi$ by replacing $e_j$ with  $e_{j+1}$ for all $1\leq j\leq s$ with the convention that $e_{s+1}=n$. 
	\item Step 2:  Let    $\alpha_{r,n}(\pi,c)=\tau=\tau_1\tau_2\ldots\tau_n$  be the permutation  obtained  from $\sigma$ by replacing $\sigma_x=\sigma_{i_y}$ with $e_1$,
	replacing $\sigma_{i_j}$ with $\sigma_{i_{j-1}}$ for all $y<j\leq k$,   and inserting $ \sigma_{i_k}$ immediately before $\sigma_{r}$.
\end{itemize}

Continuing with our running example where $\pi=836295417$, $r=5$ and $c=6$,   we have $e_1=6$, $e_2=8$,   $e_3=9$  and $\mathcal{A}(\pi)=\{1,3\}$.
Recall that the $5$-level-den-labeling    of the index sequence of  $\pi$ is given by$$
\begin{array}{lllllllllll} 
	&_{6}1&_{4}2&_{7}3&_{5}4&_{8}5&_{3}{\bold 6}&_{9}7&_{2}{\bold 8}&_{1}{\bold 9}&_{0}. 
\end{array}
$$ 
 It is easily seen that the space labeled by $6$ in the    $5$-level-den-labeling of the index sequence of $\pi$ is immediately before the index $x=1$ and $x\in \mathcal{A}(\pi)$. First we generate a permutation $\sigma=9382(10)5417$ from $\pi$ by replacing $e_1=6$ with $e_2=8$, replacing $e_2=8$ with $e_3=9$ and replacing $e_3=9$ with $e_4=10$.  Then we get the permutation $\alpha_{5,10}(\pi,6)=\tau=63928(10)5417$        by replacing $\sigma_{1}=9$ with $e_1=6$, replacing $\sigma_{3}=8$ with $\sigma_{1}=9$,   and inserting $\sigma_{3}=8$ immediately before $\sigma_{5}=10$.
\end{framed}

Now we proceed to   show that the map $\alpha_{r,n}$ verifies the following desired properties. 
\begin{lemma}\label{lemalpha}
	Let $1\leq r<n$. The map  $\alpha_{r,n}:\mathcal{X}_{n,r}\longrightarrow \mathfrak{S}^*_{n, r} $ is an injection  such that for any
	$(\pi, c)\in \mathcal{X}_{n,r} $, we have  
		\begin{equation}\label{alphare1}
		\mathrm{exc}_r(\alpha_{r,n}(\pi, c))= \mathrm{exc}_r(\pi)+1		 	 
		\end{equation}
		and 
			\begin{equation}\label{alphare2}
		\mathrm{den}_r(\alpha_{r,n}(\pi, c))=\mathrm{den}_r(\pi)+c.
		 \end{equation}
 \end{lemma}
\pf Here we use the  notations of the definition of  the map $\alpha_{r,n}$. Let $\alpha_{r,n}(\pi, c)=\tau=\tau_1\tau_2\ldots \tau_n$.
 From the construction of $\sigma$, one can easily check that $\mathrm{Excp}_r(\pi)=\mathrm{Excp}_r(\sigma)$,  $\mathrm{Nexc}_r(\pi)=\mathrm{Nexc}_r(\sigma)$, and $\mathrm{Exc}_r(\pi)$ is  order-isomorphic to $\mathrm{Exc}_r(\sigma)$. This yields that  $ 
\mathrm{den}_r(\sigma)=\mathrm{den}_r(\pi).
$

\noindent{\bf Fact 1:} For all $i\geq r$, the index $i$ is an $r$-level excedance place of $\sigma$ if and only if $i+1$ is an $r$-level excedance place of $\tau$.\\ 
By the construction of $\sigma$,  if $\sigma_i>i$ and $\sigma_i\geq r$, then we have $(i, \sigma_i)\cap \{e_1, e_{2}, \ldots, e_{s}\}\neq \emptyset$. Then  $\sigma_i>i\geq r$ implies that $\sigma_i>i+1\geq r+1$. By the construction of $\tau$,  we have $\tau_{i+1}=\sigma_i$ for all $i\geq r$. 
 Hence, we deduce that
 $\sigma_i>i\geq r$ if and only $\tau_{i+1}>i+1\geq r+1$, completing the proof of Fact 1. 

\noindent{\bf Fact 2:}  For all $i<r$, the index $i$ is an $r$-level excedance place of $\sigma$ if and only if $i$ is an $r$-level excedance place of $\tau$.\\
  Since $ \mathrm{Excp}_{r}(\pi)=\mathrm{Excp}_{r}(\sigma)$ and $\mathcal{A}(\pi)=\{i\mid \pi_i\geq r, i<r\}=\{i_1, i_2, \ldots, i_k\}$,  we have $\{i\mid \sigma_i\geq r, i<r\}=\{i_1, i_2, \ldots, i_k\}$.  Recall that $\tau$ is obtained from $\sigma$ by replacing $\sigma_x=\sigma_{i_y}$ with $e_1$,
replacing $\sigma_{i_j}$ with $\sigma_{i_{j-1}}$ for all $y<j\leq k$,   and inserting $ \sigma_{i_k}$ immediately before $\sigma_{r}$. It is easily seen that this procedure preserves the set of the $r$-level excedance places located to the left of $\sigma_r$ in $\sigma$.  This implies that $\sigma_i\geq r$ and $i<r$ if and only if $\tau_i\geq r$ and $i<r$, completing  the proof of Fact 2. 

\noindent{\bf Fact 3:} The index $r$ is an $r$-level excedance place of $\tau$. \\
  Note that  $\tau_r=\sigma_{i_k}$.
As $i_k\in \mathrm{Excp}_r(\pi)=\mathrm{Excp}_r(\sigma) $, we have 
$ 
(i_k, \sigma_{i_k})\cap \{e_1, e_2, \ldots, e_s\}\neq \emptyset 
$ by the construction of $\sigma$.
This implies that $\sigma_{i_k}>r$ and thus $\tau_r>r$. Hence, we have  $r\in \mathrm{Excp}_r(\tau)$, completing the proof of Fact 3.

Recall that 
$$
\mathrm{den}_r(\sigma)=\sum\limits_{i\in \mathrm{Excp}_r(\sigma)} i+{\mathrm{inv}}(\mathrm{Exc}_{r}(\sigma))+{\mathrm{inv}}(\mathrm{Nexc}_{r}(\sigma)). 
$$
 By Fact 1 through Fact 3,   one can easily check that the insertion of $\sigma_{i_k}$ immediately before $\sigma_r$ will increase the sum $\sum\limits_{i\in \mathrm{Excp}_r(\sigma)}i$ by $r+\mathrm{exc}_r(\sigma)$.   Moreover,  the procedure from $\sigma$ to $\tau$ does not affect $\mathrm{inv}(\mathrm{Nexc}_r(\sigma))$, and  $\mathrm{Exc}_r(\tau)$ is obtained from $\mathrm{Exc}_r(\sigma)$ by inserting $e_1$ immediately before the $y$-th letter of  $\mathrm{Exc}_r(\sigma)$ (counting from left to right). As $e_1$ is less than all the $r$-level excedance letters of $\sigma$,   the procedure from $\sigma$ to $\tau$ will increase $\mathrm{inv}(\mathrm{Exc}_r(\sigma))$ by $y-1$. 
 Recall that
 $\mathrm{den}_r(\sigma)=\mathrm{den}_r(\pi)$,  $\mathrm{exc}_r(\sigma)=\mathrm{exc}_r(\pi)$, and  $c=\mathrm{exc}_r(\pi)+r+y-1$.
 Therefore, we deduce that
 $$
 \begin{array}{ll}
 	\mathrm{den}_r(\tau)&=\sum\limits_{i\in \mathrm{Excp}_r(\tau)} i+{\mathrm{inv}}(\mathrm{Exc}_{r}(\tau))+{\mathrm{inv}}(\mathrm{Nexc}_{r}(\tau))\\
 	&=\sum\limits_{i\in \mathrm{Excp}_r(\sigma)} i+{\mathrm{inv}}(\mathrm{Exc}_{r}(\sigma))+{\mathrm{inv}}(\mathrm{Nexc}_{r}(\sigma))+\mathrm{exc}_r(\sigma)+r+y-1\\ 	  	&=\mathrm{den}_r(\sigma)+\mathrm{exc}_r(\sigma)+r+y-1\\
 	&=\mathrm{den}_r(\pi)+\mathrm{exc}_r(\pi)+r+y-1\\
 	&=\mathrm{den}_r(\pi)+c,
 \end{array}
 $$ completing the  proof of~\eqref{alphare2}.
 Again by  Fact  1 through Fact 3, the procedure from $\sigma$ to $\tau$ will increase $\mathrm{exc}_r(\sigma)$ by $1$. As $\mathrm{exc}_r(\pi)=\mathrm{exc}_r(\sigma)$, we have
 $ 
 \mathrm{exc}_r(\tau)=\mathrm{exc}_r(\pi)+1,
 $ 
  completing the  proof of~\eqref{alphare1}.

   Next we proceed to  show  that $\tau\in \mathfrak{S}^*_{n, r} $. 
      
   \noindent{\bf Fact 4.} $\mathrm{Excl}_r( \tau)=\{e_1, e_2, \ldots, e_{s+1}\}$.\\
   By  Fact  1 through Fact 3, it is not difficult to check that the procedure from $\sigma$ to $\tau$ will generate a new $r$-level excedance letter  $e_1$   and preserve all the $r$-level excedance letters of $\sigma$. By the definition of the procedure from $\pi$ to $\sigma$, it is easily checked  that $\mathrm{Excl}_r( \sigma)=\{ e_2,e_3, \ldots, e_{s+1}\}$. Therefore,  we deduce that $\mathrm{Excl}_r( \tau)=\{e_1, e_2, \ldots, e_{s+1}\}$ as desired, completing the proof of Fact~4.
   
   By Fact~3, we have $\tau_r>r$. Hence, in order to prove that $\tau\in \mathfrak{S}^*_{n, r} $, it remains to show that the index $i$ is not an $r$-level grande-fixed place for all $i\geq r$.       
  Take $t>r$ arbitrarily (if any) such that $\tau_t\geq t$.   
  By Fact 1, $\tau_t= t$ would imply  that $\sigma_{t-1}\leq t-1$.  However, 
  by the construction of $\tau$,  we have $ \sigma_{t-1}=\tau_t=t$,  a contradiction. Hence, we have $\tau_t>t$. Again by Fact~1, $\sigma_{t-1}$ is an $r$-level excedance letter in $\sigma$.   By the construction of $\sigma$, one can easily check that 
  $ 
  (t-1, \sigma_{t-1})\cap \{e_1, e_{2}, \ldots, e_s\}\neq \emptyset 
  $.
  Then   we have   $[t, \tau_{t})\cap \mathrm{Excl}_r(\tau)=(t-1, \sigma_{t-1})\cap \mathrm{Excl}_r(\tau)\neq \emptyset $ as $\mathrm{Excl}_r( \tau)=\{e_1, e_2, \ldots, e_{s+1}\}$ by Fact~4.   This proves that  for all $i>r$, the index $i$ is not an $r$-level grande-fixed place of $\tau$. 
  As $\tau_r>r$,  $e_1\in \mathrm{Excl}_r( \tau)$ and $\tau_x=e_1$, we have
  $e_1\in [r, \tau_r)\cap \mathrm{Excl}_r( \tau)$. This implies that
  the index $r$ is not  an $r$-level grande-fixed place of $\tau$. Hence, we conclude that $\tau\in \mathfrak{S}^*_{n, r} $.

  In the following, we aim to show that the map $\alpha_{r,n}$ is an injection. By (\ref{alphare2}), we have $c=\mathrm{den}_r(\tau)-\mathrm{den}_r(\pi)$. Therefore,  in order to  show that $\alpha_{r,n}$ is injective, it suffices to show that we can recover $\pi$ from $\tau$. 
    First we shall prove that the resulting permutation    $\tau$ verifies the following desired  properties.
\begin{itemize}
	\item[\upshape(\rmnum{1})]   $\tau_x=e_1$ is the smallest  $r$-level excedance letter of $\tau$.
	\item[\upshape(\rmnum{2})]  The set of the $r$-level excedance letters that occur weakly  to the right of $\tau_x$ and  weakly to the left of $\tau_r$  is given by 
	$\{\tau_{i_j}\mid y\leq j\leq k+1\}$ with $x=i_{y}<i_{y+1}<\cdots< i_k<i_{k+1}=r$  and $\tau_{i_j}=\sigma_{i_{j-1}}$ for all $y<j\leq k+1$. 
	\item[\upshape(\rmnum{3})] $\mathrm{Excl}_r(\tau)-\mathrm{Excl}_r(\sigma)=\{\tau_x\}$ and $\mathrm{Excl}_r(\tau)=\{e_1, e_2, \ldots, e_{s+1}\}$ with $e_{s+1}=n$. 
\end{itemize} 
  Clearly,  $e_1=\tau_x$ is the smallest $r$-level excedance letter of $\tau$ since  $\tau_x=e_1\in \mathrm{Excl}_r( \tau)=\{e_1, e_2, \ldots, e_{s+1}\}$ with $e_1<e_2<\ldots<e_{s+1}$, completing the proof of (\rmnum{1}).  As $\mathrm{Excp}_r(\pi)=\mathrm{Excp}_r(\sigma)$, we have
 $\{j\mid \sigma_i\geq r, i<r\}=\{i\mid \pi_i\geq r, i<r\}=\mathcal{A}(\pi)=\{i_j\mid 1\leq j\leq k\}$. Then   Fact 2 and Fact 3 ensure  that 
 the set of the $r$-level excedance letters that occur weakly to the right of $\tau_x$ and  weakly to the left of $\tau_r$  is given by 
 $\{\tau_{i_j}\mid y\leq j\leq k+1\}$ where $x=i_{y}<i_{y+1}<\ldots< i_k<i_{k+1}=r$. Then by the definition of the procedure from $\sigma$ to $\tau$,  we have $\tau_{i_j}=\sigma_{i_{j-1}}$ for all $y<j\leq k+1$,  completing the proof of  (\rmnum{2}). The second part of (\rmnum{3}) follows directly from Fact 4. According to the definition of  the procedure from $\pi$ to $\sigma$, it is easy to check that $\mathrm{Excl}_r( \sigma)=\{e_1, e_2, \ldots, e_{s+1}\}-\{e_1\}$, completing the proof of (\rmnum{3}).  
 
 Now we are in the position to give  a description of the procedure  to recover the permutation $\pi$ from $\tau$.
 \begin{itemize}
 	\item  Choose $\tau_{w}$ to be  the smallest  $r$-level excedance letter of $\tau$. 
 	\item 
 	Find all the $r$-level excedance letters,   say   $\tau_{j_1}, \tau_{j_2}, \ldots, \tau_{j_\ell}$ with $w=j_1<j_2<\ldots<j_{\ell}=r$,  that are  located  weakly  to the right of $\tau_{w}$ and weakly to the left of  $\tau_r$ in $\tau$.  
 	Construct a permutation $\sigma$ from $\tau$ by removing $\tau_{j_\ell}$ and then replacing  each $\tau_{j_{i-1}}$  with   $\tau_{j_{i}}$ for all $1< i\leq \ell$.
 	\item  Find all the $r$-level excedance  letters of $\sigma$, say $e'_1, e'_2, \ldots, e'_t$ with $e'_1<e'_2<\ldots<e'_t$.
   Let $\pi$ be the permutation obtained from $\sigma$ by replacing each $e'_i$ with $e'_{i-1}$ for all $1\leq i\leq t$ with the convention that $e'_0=\tau_w$.
 \end{itemize} 
 Properties (\rmnum{1})-(\rmnum{3}) ensure  that after applying  the above procedure to $\tau$, we  will result in the same permutation $\sigma$ and thus the same permutation $\pi$. This completes the proof.\qed

\begin{framed}
\begin{center}
{\bf The map $\beta_{r,n}:\mathcal{Y}_{n,r}\longrightarrow \mathfrak{S}_{n, r}-\mathfrak{S}^*_{n, r} $} 
\end{center}
 Let $(\pi, c)\in \mathcal{Y}_{n,r}$ with   $|\mathrm{Exc}_r(\pi)|=s$ and let  $0<c\leq n-1$.  Assume that   the space labeled by $c$  under the $r$-level-den-labeling of the index sequence of $\pi$    is  immediately before the index $x$   and that $\mathrm{Excl}_r(\pi)=\{e_1,e_2, \ldots, e_{s}\}$ with $e_1<e_2\ldots<e_s$.    
 Let  $\mathcal{B}(\pi)=\{j_1, j_2, \ldots, j_{\ell}\}$
 and that $x=j_y$ for some $1\leq y\leq \ell$. 
 Define  $\beta_{r,n}(\pi,c)$ to be the permutation  constructed by the following procedure.
 \begin{itemize}
 	\item Step 1: Construct a permutation  $\sigma=\sigma_1\sigma_2\ldots \sigma_{n-1}$  from $\pi$ by replacing $e_j$ with  $e_{j+1}$ for all $1\leq j\leq s$ with the convention that $e_{s+1}=n$. 
 	\item Step 2: Construct a permutation  
 	 $\sigma'=\sigma'_1\sigma'_2\ldots \sigma'_{n}$   from $\sigma$ by    replacing $\sigma_x=\sigma_{j_y}$ with $e_1$,
 	replacing $\sigma_{j_i}$ with $\sigma_{j_{i-1}}$ for all $y<i\leq \ell$,   and inserting $ \sigma_{j_\ell}$ immediately before $\sigma_{r}$.	
 
 	\item  Step 3: Find all the $r$-level excedance letters of $\sigma'$ that occur to the left of $\sigma'_r$, say $e'_1, e'_2, \ldots, e'_k$ with  $e'_1<e'_2<\ldots< e'_k$.  
 	Then define   $ \beta_{r,n}(\pi,c)$ to be the permutation obtained from $\sigma'$ by replacing
 	$e'_j$ with $e'_{j-1}$ for all $1< j\leq k$ and replacing $e'_1$ with $e'_k$. 
 \end{itemize}

 Continuing with our running example where $\pi=836295417$, $r=5$ and $c=4$, we have $e_1=6$, $e_2=8$,  $e_3=9$ and $\mathcal{B}(\pi)=\{2,4\}$. Recall that the $5$-level-den-labeling    of the index sequence of  $\pi$ is given by$$
 \begin{array}{lllllllllll} 
 	&_{6}1&_{4}2&_{7}3&_{5}4&_{8}5&_{3}{\bold 6}&_{9}7&_{2}{\bold 8}&_{1}{\bold 9}&_{0}. 
 \end{array}
 $$ It is easily seen that the space labeled by $4$ in the    $5$-level-den-labeling of the index sequence of $\pi$ is immediately before the index $x=2$ and $x\in \mathcal{B}(\pi)$. First we generate a permutation $\sigma=9382(10)5417$  from $\pi$ by replacing $e_1=6$ with $e_2=8$, replacing $e_2=8$ with $e_3=9$ and replacing $e_3=9$ with $e_4=10$.   Then we generate a permutation $\sigma'=96832(10)5417$ from $\sigma$ by replacing $\sigma_2=3$ with $e_1=6$ and replacing $\sigma_4=2$ with $\sigma_2=3$ and inserting $\sigma_4=2$ immediately before $\sigma_5=10$. Clearly, all the $5$-level excedance letters that occur to the left of $\sigma'_5=2$ in $\sigma'$ are given by $e'_1=6, e'_2=8, e'_3=9$.
 Finally, we get the permutation  $\beta_{5,10}(\pi,4)= 89632(10)5417$   obtained from $\sigma'$ by   replacing $e'_1=6$ with $e'_3=9$,  replacing $e'_2=8$ with $e'_1=6$ and replacing $e'_3=9$ with $e'_2=8$. 
 
 \end{framed}

\begin{lemma}\label{lembeta}
	Let $1\leq r<n$. The map  $\beta_{r,n}:\mathcal{Y}_{n,r}\longrightarrow \mathfrak{S}_{n, r}-\mathfrak{S}^*_{n, r} $ is an injection  such that for any
	$(\pi, c)\in \mathcal{Y}_{n,r} $, we have  
	\begin{equation}\label{betare1}
		\mathrm{exc}_r(\beta_{r,n}(\pi, c))= \mathrm{exc}_r(\pi)		 	 
	\end{equation}
	and 
	\begin{equation}\label{betare2}
		\mathrm{den}_r(\beta_{r,n}(\pi, c))=\mathrm{den}_r(\pi)+c.
	\end{equation}
\end{lemma}
 \pf  Here we use the notations of   the definition of  the map $\beta_{r,n}$. Let $\beta_{r,n}(\pi, c)=\tau=\tau_1\tau_2\ldots \tau_n$. By the construction of $\sigma$, it is easily seen that we have $\mathrm{Excp}_r(\pi)=\mathrm{Excp}_r(\sigma)$,  $\mathrm{Nexc}_r(\pi)=\mathrm{Nexc}_r(\sigma)$, and $\mathrm{Exc}_r(\pi)$ is  order-isomorphic to $\mathrm{Exc}_r(\sigma)$. This ensures that  $\mathrm{den}_r(\sigma)=\mathrm{den}_r(\pi)$. 
  
 \noindent{\bf Fact ${\bf1'}$.} The index $i$ is an $r$-level excedance place of $\sigma$ if and only if $i+1$ is an $r$-level excedance place of $\sigma'$ for all $i\geq r$.\\
 One can  verify Fact $1'$ by the same reasoning as in the proof of Fact 1 and the proof is omitted here.
 
 By the definition of the procedure from $\pi$ to $\sigma$, it is easily seen that  $\mathrm{Excl}_r(\sigma)=\{e_1, e_2, \ldots, e_{s+1}\}-\{e_1\}$. Since  $\mathrm{Excp}_r(\pi)=\mathrm{Excp}_r(\sigma)$, we have 
  $$\{j\mid \sigma_j<r, j<r\}=\{i\mid \pi_i<r, i<r\}=\mathcal{B}(\pi)=\{j_1, j_2, \ldots, j_\ell\}.$$
  Namely, we have  $\sigma_{j_i}<r$ for all $1\leq i\leq \ell$.  Then by the definition of   the procedure from $\sigma$ to $\sigma'$, we deduce the following three facts.
  
 \noindent{\bf Fact ${\bf2'}$.} For all $i<r$ and $i\neq x$, the index $i$ is an $r$-level excedance place of $\sigma$ if and only if  $i$  is an $r$-level excedance place of $\sigma'$.\\  
  \noindent{\bf Fact ${\bf3'}$.} The index $x$ is an $r$-level excedance place of $\sigma'$ and the index  $r$ is not an $r$-level excedance place of $\sigma'$.\\
   \noindent{\bf Fact ${\bf4'}$.} $\mathrm{Excl}_r(\sigma')=\{e_1, e_2, \ldots, e_{s+1}\}$.

  Recall that 
  $$
  \mathrm{den}_r(\sigma)=\sum\limits_{i\in \mathrm{Excp}_r(\sigma)} i+{\mathrm{inv}}(\mathrm{Exc}_{r}(\sigma))+{\mathrm{inv}}(\mathrm{Nexc}_{r}(\sigma)). 
  $$
  By  Fact $1'$ through Fact $3'$,   one can easily check that the insertion of $\sigma_{j_\ell}$ immediately before $\sigma_r$ and replacing $\sigma_x$ with $e_1$ will increase the sum $\sum\limits_{i\in \mathrm{Excp}_r(\sigma)}i$ by $x+\mathrm{exc}_r(\sigma)$.  Moreover,   the procedure from $\sigma$ to $\sigma'$ does not affect $\mathrm{inv}(\mathrm{Nexc}_r(\sigma))$.

   Assume that there are $p$ $r$-level excedance letters located  to the left of $\pi_x$ and $q$ $r$-level excedance letters located  to the right of $\pi_x$ and to the left of $\pi_r$ in $\pi$. By the rules specified  in the $r$-level-den-labeling of the places of the index sequence of $\pi$, we have $c=\mathrm{exc}_r(\pi)+p+q+y$ and $x=p+y$. Since $\mathrm{Excp}_r(\sigma)=\mathrm{Excp}_r(\pi)$,  we deduce that there are $p$ $r$-level excedance letters located  to the left of $\sigma_x$ and  $q$ $r$-level excedance letters located  to the right of $\sigma_x$ and to the left of $\sigma_r$ in $\pi$.  Then replacing $\sigma_x$ with $e_1$ in the procedure from $\sigma$ to $\sigma'$ will increase $\mathrm{inv}(\mathrm{Exc}_r(\sigma))$ by $p$  since $e_1$ is less than all the $r$-level  excedance letters of $\sigma$.  
   Therefore, we derive  that
   \begin{equation}\label{eqbeta1}
   	\begin{array}{ll}
   		\mathrm{den}_r(\sigma')&=\sum\limits_{i\in \mathrm{Excp}_r(\sigma')} i+{\mathrm{inv}}(\mathrm{Exc}_{r}(\sigma'))+{\mathrm{inv}}(\mathrm{Nexc}_{r}(\sigma'))\\
   		&=\sum\limits_{i\in \mathrm{Excp}_r(\sigma)} i+{\mathrm{inv}}(\mathrm{Exc}_{r}(\sigma))+{\mathrm{inv}}(\mathrm{Nexc}_{r}(\sigma))+\mathrm{exc}_r(\sigma)+x+p\\
   		&=	\mathrm{den}_r(\sigma)+\mathrm{exc}_r(\sigma)+p+x\\
   		&=	\mathrm{den}_r(\pi)+\mathrm{exc}_r(\pi)+p+x,\\
   	\end{array}
   \end{equation}
   where the last equality follows from the fact that $\mathrm{exc}_r(\pi)=\mathrm{exc}_r(\sigma)$ and $\mathrm{den}_r(\pi)=\mathrm{den}_r(\sigma)$.

 It is apparent that the procedure from $\sigma'$ to $\tau$  does not affect $\mathrm{inv}(\mathrm{Nexc}(\sigma'))$ and the sum $\sum\limits_{i\in \mathrm{Excp}_r(\sigma')}i$. 
 Recall that $\sigma'_x=e_1$ is the smallest $r$-level excedance letter of $\sigma'$.  
 This implies that $e'_1=\sigma'_x$. 
 Recall that there are $p$ $r$-level excedance letters located  to the left of $\sigma_x$ and  $q$ $r$-level excedance letters located  to the right of $\sigma_x$ and to the left of $\sigma_r$ in $\pi$.  Then by Facts $2'$ and $3'$, the number of $r$-level excedance letters located  to the left of $\sigma'_x$  in $\sigma'$ is given by $p$ and the number of 
 $r$-level excedance letters located  to the right of $\sigma'_x$ and to the left of $\sigma'_r$ in $\sigma'$ is given by $q$. 
 Thus, by the definition of the procedure from $\sigma'$ to $\tau$, it is not difficult to check  that  
 \begin{equation}\label{eqbeta3}
 	\mathrm{inv}(\mathrm{Exc}(\tau))=\mathrm{inv}(\mathrm{Exc}(\sigma'))-p+q.
 \end{equation}
 
  Recall that $\mathrm{exc}_r(\sigma)=\mathrm{exc}_r(\pi)$, $\mathrm{den}_r(\sigma)=\mathrm{den}_r(\pi)$, $x=p+y$, and   $c=\mathrm{exc}_r(\pi)+p+y+q$.
 Invoking (\ref{eqbeta1}) and (\ref{eqbeta3}), we derive that
 $$
 \begin{array}{ll}
 	\mathrm{den}_r(\tau)&=\sum\limits_{i\in \mathrm{Excp}_r(\tau)} i+{\mathrm{inv}}(\mathrm{Exc}_{r}(\tau))+{\mathrm{inv}}(\mathrm{Nexc}_{r}(\tau))\\
 	&=\sum\limits_{i\in \mathrm{Excp}_r(\sigma')} i+{\mathrm{inv}}(\mathrm{Exc}_{r}(\sigma'))+{\mathrm{inv}}(\mathrm{Nexc}_{r}(\sigma')-p+q\\ 
 	&=\mathrm{den}_r(\sigma')-p+q\\
 	  	&= \mathrm{den}_r(\pi)+\mathrm{exc}_r(\pi)+p+y+q\\
 	&= \mathrm{den}_r(\pi)+c,
 \end{array}
 $$
 completing the proof of  (\ref{betare2}).

 Clearly, the procedure from $\sigma'$ to $\tau$ preserves the number of the $r$-level excedance letters that occur weakly to the right of $\sigma'_r$. By Fact $1'$ through Fact $3'$, the procedure from $\sigma$ to $\sigma'$ also preserves the number of the $r$-level excedance letters that occur weakly to the right of $\sigma_r$. Hence, we derive that 
  $$
  \mathrm{exc}_r(\tau) = \mathrm{exc}_r(\sigma')=\mathrm{exc}_r(\sigma)= \mathrm{exc}_r(\pi), 
 $$
 where the last equality follows from  the fact that $\mathrm{exc}_r(\sigma)=\mathrm{exc}_r(\pi)$, 
 completing the proof of (\ref{betare1}).

  Next we proceed to  prove that $\tau\in \mathfrak{S}_{n, r}- \mathfrak{S}^*_{n, r} $.   By the definition of the procedure from $\sigma'$ to $\tau$, we have $\sigma'_i=\tau_i$ for all $i\geq r$. Moreover,  the procedure from $\sigma'$ to $\tau$ preserves the set of $r$-level  excedance letters and hence preserves the set of  $r$-level grande-fixed places.  By Fact $3'$, we have $\sigma'_r<r$. Hence, in order to prove that $\tau\in  \mathfrak{S}_{n, r}- \mathfrak{S}^*_{n, r}$, it remains to show that the index $i$ is not an $r$-level grande-fixed place of $\sigma'$ for any $i>r$.        
 Take $t>r$ arbitrarily (if any) such that $\sigma'_t\geq t$. By the same reasoning as in the proof of Lemma~\ref{lemalpha},  one can deduce that $\sigma'_t>t$ by Fact $1'$. Again by Fact $1'$, $\sigma_{t-1}$ is an $r$-level excedance letter in $\sigma$.   By the definition of the procedure from $\pi$ to $\sigma$, one can easily check that 
 $ 
 (t-1, \sigma_{t-1})\cap \{e_1, e_{2}, \ldots, e_s\}\neq \emptyset 
 $.
 By Fact $4'$, we have $\mathrm{Excl}_r( \sigma')=\{e_1, e_2, \ldots, e_{s+1}\}$. 
 Hence, we have  $[t, \sigma'_{t})\cap \mathrm{Excl}_r(\sigma')=(t-1, \sigma_{t-1})\cap \mathrm{Excl}_r(\sigma')\neq \emptyset $, completing the proof of the fact that $\tau\in \mathfrak{S}_{n, r}- \mathfrak{S}^*_{n, r}$.

 In the following, we aim to show that the map $\beta_{r,n}$ is injective. By (\ref{betare2}), we have $c=\mathrm{den}_r(\tau)-\mathrm{den}_r(\pi)$. In order to  show that $\beta_{r,n}$ is injective, it is sufficient  to show that we can recover $\pi$ from $\tau$. 
 First we shall prove that the resulting permutations    $\tau$ and $\sigma'$ verify the following desired  properties.
 \begin{itemize}
 	\item[\upshape(\rmnum{1}$^{'}$)]   $\tau_x$ is the greatest   $r$-level excedance letter that occurs to the left of $\tau_r$ in $\tau$.
 	\item[\upshape(\rmnum{2}$^{'}$)]   The smallest  $r$-level excedance letter  of $\tau$ occurs to the left of $\tau_r$ in $\tau$. 	 	
 	 		\item[\upshape(\rmnum{3}$^{'}$)]   $\sigma'_x=e_1$ is the smallest  $r$-level excedance letter of $\sigma'$. 
 	\item[\upshape(\rmnum{4}$^{'}$)]  The set of the non-$r$-level excedance letters that occur   to the right of $\sigma'_x$ and  weakly to the left of $\sigma'_r$  is given by 
 	$\{\sigma'_{j_i}\mid y<i\leq \ell+1\}$ with $x<i_{y+1}<\ldots< i_{\ell+1}=r$  and $\sigma'_{j_i}=\sigma_{j_{i-1}}$ for all $y<i\leq \ell+1$. 
 	\item[\upshape(\rmnum{5}$^{'}$)] $\mathrm{Excl}_r(\sigma')-\mathrm{Excl}_r(\sigma)=\{\sigma'_x\}$ and $\mathrm{Excl}_r(\sigma')=\{e_1, e_2, \ldots, e_{s+1}\}$ with $e_{s+1}=n$. 
 \end{itemize} 
  Relying on Facts $1'$-$4'$, one can prove 
  (\rmnum{3}$^{'}$)-(\rmnum{5}$^{'}$) by similar arguments as in the proof of (\rmnum{1})-(\rmnum{3}) of Lemma \ref{lemalpha}. 
  By (\rmnum{3}$^{'}$), we have $e'_1=\sigma'_x=e_1$. Then (\rmnum{1}$^{'}$) and (\rmnum{2}$^{'}$) follow directly from the definition of the procedure from $\sigma'$ to $\tau$.  
 
 Now we are in the position to give  a description of the procedure  to recover the permutation $\pi$ from $\tau$.
 \begin{itemize}
 	\item  Choose $\tau_{w}$ to be  the greatest   $r$-level excedance letter of $\tau$ that occurs to the left of $\tau_r$. 
 	\item  Find all the  $r$-level excedance letters of $\tau$ that occurs to the left of $\tau_r$, say $f_1, f_2, \ldots, f_t $ with $f_1<f_2<\cdots<f_t=\tau_w$. 	Construct a permutation $\sigma'$ by replacing $f_{i-1}$ with $f_{i}$ for all $1<i\leq t$ and replacing $f_t$ with $f_1$. 
 	\item  
 	Find all the non-$r$-level excedance letters,   say   $\sigma'_{i_1}, \sigma'_{i_2}, \ldots, \sigma'_{i_m}$ with $w<i_1<i_2<\cdots<i_{m}=r$,  that are  located  to the right of $\sigma'_{w}$ and weakly to the left of  $\sigma'_r$ in $\sigma'$.  
 	Construct a permutation $\sigma$ by removing $\sigma'_{i_m}$ and then replacing  each $\sigma'_{i_{j-1}}$  with   $\sigma'_{i_{j}}$ for all $1\leq j\leq m$ with the convention that $i_0=w$. 
 	\item  Find all the $r$-level excedance  letters of $\sigma$, say $g_1, g_2, \ldots, g_h$ with $g_1<g_2<\cdots<g_h$.
 	Let $\pi$ be the permutation obtained from $\sigma$ by replacing  $g_i$ with $g_{i-1}$ for each $1\leq i\leq h$ with the convention that $g_0=\sigma'_w$.
 \end{itemize} 
 Properties (\rmnum{1}$^{'}$)-(\rmnum{5}$^{'}$) ensure that after applying  the above procedure to $\tau$, we  will result in the same permutation $\sigma$,  the same $\sigma'$ and thus the same permutation $\pi$. This completes the proof.\qed

   Let $\pi\in \mathfrak{S}_{n-1}$ with   $|\mathrm{Exc}_r(\pi)|=s$ and $0\leq c\leq n-1$. Assume that   the space labeled by $c$   under the $r$-level-den-labeling of the index sequence of $\pi$    is  immediately before the index $x$  for some $ x\geq r$ (when $c >0$)  and that $\mathrm{Excl}_r(\pi)=\{e_1,e_2, \ldots, e_{s}\}$ with $e_1<e_2\ldots<e_s$.  Let $e_{s+1}=n$. 
     Remark~\ref{ob1}  tells us that the space before the index $x$   is also labeled by $c$ when $c>0$  and the   space after $n-1$ is also labeled by $0$  under  the den-labeling of the index sequence of $\pi$.
    In the following, we proceed to show that the resulting permutation $\tau=\phi_n(\pi, c)$ has the following properties.

  \begin{lemma}\label{lemtheta}
  	For $c>0$ and $x\geq r$ or $c=0$, the resulting permutation $\tau=\phi_n(\pi, c)$   verifies the following properties. 
  	\begin{itemize}
  		\item [\upshape(\rmnum{1})]  We have 
  		$$ 
  		\mathrm{exc}_r(\tau)=\left\{ \begin{array}{ll}
  			\mathrm{exc}_r(\pi)&\, \mathrm{if}\,\, 0\leq c\leq \mathrm{exc}_r(\pi)+r-1,\\
  			\mathrm{exc}_r(\pi)+1&\, \mathrm{otherwise}, 
  		\end{array}
  		\right.
  		$$
  		and 
  		$$
  		\mathrm{den}_r(\tau)=\mathrm{den}_r(\pi)+c.
  		$$
  		\item[\upshape(\rmnum{2})]    We have $\tau\in \mathfrak{S}_{n}-\mathfrak{S}_{n,r}$. 
  		
  	\end{itemize} 
  	
  \end{lemma}
\pf It is apparent that the statement holds for $c=0$. Now we assume that $c>0$. As $x\geq r$, it follows from Lemma~\ref{Han}~(ii) that $x$ is an $r$-level grande-fixed place (greatest one) of $\tau$, which proves~(\rmnum{2}). 
In the following, we proceed to verify (\rmnum{1}).

Suppose that the insertion letter for the  space before the index $x$ is $e_{y}$, where $1\leq y\leq s+1$. Recall that the permutation $\tau=\phi_n(\pi, c)$ is constructed by the following procedure.
\begin{itemize}
	\item  First construct a permutation $\sigma=\sigma_1\sigma_2\ldots \sigma_{n-1}$ by replacing $e_j$ with  $e_{j+1}$ for all $y\leq j\leq s $.
	\item Then $\tau$ is obtained    obtained from $\sigma$ by inserting the letter $e_y$ immediately before the letter $\sigma_{x}$.  
\end{itemize}
  From the construction of $\sigma$, one can easily check that $\mathrm{Excp}_r(\pi)=\mathrm{Excp}_r(\sigma)$ ,  $\mathrm{Nexc}_r(\pi)=\mathrm{Nexc}_r(\sigma)$, and $\mathrm{Exc}_r(\pi)$ is  order-isomorphic to $\mathrm{Exc}_r(\sigma)$. This guarantees that 
   $\mathrm{den}_r(\sigma)=\mathrm{den}_r(\pi)$.
  Similarly, we have 
  $\mathrm{den}(\sigma)=\mathrm{den}(\pi)$.
  
   \noindent{\bf Fact A.} For each $i\geq x$, the index $i$ is an  excedance  place (resp., an  $r$-level excedance  place) of $\sigma$ if and only if the index  $i+1$ is  an  excedance  place (resp., an  $r$-level excedance  place) of $\tau$; for each $i<x$, the index  $i$ is  an  excedance  place (resp., an  $r$-level excedance  place) of $\sigma$ if and only if the index $i$ is  an  excedance  place (resp., an  $r$-level excedance  place)  of $\tau$.
   
  One could verify the first part of  Fact A by the same reasoning as in the proof of Fact 1 and the proof is omitted here. The second part of  Fact A follows directly from the fact that $\sigma_i=\tau_i$ for all $i<x$.

  By the construction of $\sigma$,   we have $\tau_x=x$ when $x=e_y$, and $\tau_x=e_y>x$, otherwise.  Hence, the following fact  holds.

  \noindent{\bf Fact B.} The index 
  $x $ is  an  excedance  place (resp., an  $r$-level excedance  place) of $\tau$ if and only if $x\neq e_y$.

    Recall that $\mathrm{den}_r(\sigma)=
  \mathrm{den}_r(\pi)$  and  $\mathrm{den}(\sigma)=
  \mathrm{den}(\pi)$. 
  Combining Fact A and  Fact B,    it is not difficult to check that
  $$
  \mathrm{den}_r(\tau)- \mathrm{den}_r(\pi)=\mathrm{den}_r(\tau)- \mathrm{den}_r(\sigma)=\mathrm{den}(\tau)- \mathrm{den}(\sigma)=\mathrm{den}(\tau)- \mathrm{den}(\pi).
  $$   
  Then by Lemma \ref{Han}~(i), we deduce that
  $$
  \mathrm{den}_r(\tau)=\mathrm{den}_r(\pi)+c$$ as desired.
  
Recall that $\tau$ is obtained from $\sigma$ by inserting $e_y$  immediately before $\sigma_x$. By Facts A and B, it is easily seen that the procedure from $\sigma$ to $\tau$ will increase the number of $r$-level excedance places by one when $x\neq e_y$, and will preserve the number of  $r$-level excedance places, otherwise.   According to the rules specified in the $r$-level-den-labeling of the index sequence of $\pi$, we have  $c\leq \mathrm{exc}_r(\pi)+r-1$ when $x=e_y$, and  $c> \mathrm{exc}_r(\pi)+r-1$, otherwise. Then, we deduce that
$$  
	\mathrm{exc}_r(\tau)=\left\{ \begin{array}{ll}
		\mathrm{exc}_r(\sigma)=\mathrm{exc}_r(\pi)&\, \mathrm{if}\,\, 0\leq c\leq \mathrm{exc}_r(\pi)+r-1,\\
		&\\
		\mathrm{exc}_r(\sigma)+1=\mathrm{exc}_r(\pi)+1&\, \mathrm{otherwise}, 
	\end{array}
	\right. 
$$
where the last equality follows from the fact that $\mathrm{exc}_r(\pi)=\mathrm{exc}_r(\sigma)$. This completes the proof. \qed

  Now we are in the position to complete the proof of  Theorem \ref{thden}.

 \noindent{\bf Proof of Theorem \ref{thden}}.  
 Let $\pi\in \mathfrak{S}_{n-1}$  and $0\leq c\leq n-1$.
 Assume that the space labeled by $c$ in the $r$-level-den-labeling of the index sequence of $\pi$ is immediately before the index $x$ when $c\neq 0$. Then define
  $$  
  \phi^{\mathrm{den}}_{r,n}(\pi,c)=\left\{ \begin{array}{ll}    
  \alpha_{r,n}(\pi, c)&\, \mathrm{if}\,\, c>0 \,\, \mbox{and}\,\, x\in \mathcal{A}(\pi),\\
   \beta_{r,n}(\pi, c)&\, \mathrm{if}\,\, c>0 \,\, \mbox{and}\,\, x\in \mathcal{B}(\pi),\\
    \phi_n(\pi, c)&\, \mathrm{if}\,\, c>0 \,\, \mbox{and}\,\,x\geq r\,\, \mbox{or}\,\, c=0. 
  \end{array}
  \right.
 $$
  By Lemmas~\ref{lemalpha}-\ref{lemtheta}, the map   $\phi^{\mathrm{den}}_{r,n}$ maps the  pair $(\pi, c)$ to a permutation $\tau\in \mathfrak{S}_{n}$ such that the resulting permutation $\tau$ verifies   (\ref{eqprop1}) and (\ref{eqprop2}).  
By cardinality reasons,
in order to prove that $\phi^{\mathrm{den}}_{r,n}$ is a bijection,
it is sufficient to prove that $\phi^{\mathrm{den}}_{r,n}$ is injective.  This can be justified by Lemma \ref{Han} and  Lemmas~\ref{lemalpha}-\ref{lemtheta}. 
 \qed

   Notice that when $r\geq n$, we have
 $\mathrm{rdes}(\pi)=0=\mathrm{exc}_r(\pi)$ and 
  $\mathrm{rmaj}(\pi)=\mathrm{inv}(\pi)=\mathrm{den}_r(\pi)$   for any $\pi\in  \mathfrak{S}_{n} $. Therefore,   $(\mathrm{rdes}, \mathrm{rmaj})$ and $(\mathrm{exc}_r, \mathrm{den}_r)$ are equidistributed over $\mathfrak{S}_{n}$ when $r\geq n$. For $1\leq r<n$, combining   Lemma \ref{Raw} and Theorem \ref{thden}, we are led to a bijective   proof of Conjecture~\ref{con1}.

\section*{Acknowledgments}
The work  was supported by
the National Natural
Science Foundation of China grants 12071440, 12322115 and 12271301.

\end{document}